\documentclass[11pt]{article}
\usepackage{epic,latexsym,amssymb}
\usepackage{tikz}
\usepackage[bold,full]{complexity}
\usepackage{tikz,epsf}
\usepackage{hyperref}
\usepackage{caption,subcaption}
\usepackage{tabularx}
\usepackage{array}

\usepackage{amsfonts}
\usepackage{amscd}
\usepackage{amsmath}
\usepackage{graphicx}
\usepackage{color}
\usepackage{caption,subcaption}
\usepackage{lineno}
\usetikzlibrary{snakes}

\newtheorem{staticthm}{Theorem}
\newcounter{staticthmctr}
\newtheorem{thm}{Theorem}
\newtheorem{conj}[thm]{Conjecture}

\newtheorem{prop}[thm]{Proposition}
\newtheorem{lem}[thm]{Lemma}
\newtheorem{quest}{Question}
\newtheorem{problem}{Problem}

\newtheorem{cor}[thm]{Corollary}
\newtheorem{remark}[thm]{Remark}

\newtheorem{theorem}{Theorem}

\newtheorem{corollary}{Corollary}

\newcommand{\qed}{$\Box$}

\newenvironment{proof}{\noindent\textbf{Proof. }}{\hfill $\square$}

\newcommand{\vertex}{\node[vertex]}
\tikzstyle{vertex}=[circle, draw, inner sep=0pt, minimum size=6pt]

\newcommand{\QEDmark}{\mbox{\textsc{qed}}}
\newcommand{\proofStarter}[1]{\textsc{#1} }

\def\vertex(#1){\put(#1){\circle*{2}}}
\def\vertexo(#1){\put(#1){\circle{2}}}
\def\vert(#1){\put(#1){\circle*{1.5}}}
\def\verto(#1){\put(#1){\circle{1.5}}}
\def\lab(#1)#2{\put(#1){\makebox(0,0)[c]{#2}}}
\definecolor{DarkGreen}{rgb}{0.2, 0.6, 0.3}

\definecolor{electricindigo}{rgb}{0.44, 0.0, 1.0}

\begin{document}

\title{Comparing the $p$-independence number of regular graphs to the $q$-independence number of their line graphs}
\author{
Yair Caro \\
Department of Mathematics \\
University of Haifa--Oranim\\
Tivon 36006, Israel \\
\small {\tt Email: yacaro@kvgeva.org.il}\\
\\
Randy Davila\\
Department of Computational Applied \\ Mathematics \& Operations Research\\
Rice University \\
Houston, TX 77005, USA \\
\small {\tt Email: rrd6@rice.edu}\\
\\
Ryan Pepper\\
Department of Mathematics and Statistics \\
University of Houston--Downtown \\
Houston, TX 77002, USA \\
\small {\tt Email: pepperr@uhd.edu}\\
}

\date{}
\maketitle

\begin{abstract} 
Let $G$ be a simple graph and let $L(G)$ denote the \emph{line graph} of $G$. A \emph{$p$-independent} set in $G$ is a set of vertices $S \subseteq V(G)$ such that the subgraph induced by $S$ has maximum degree at most $p$. The \emph{$p$-independence number} of $G$, denoted by $\alpha_p(G)$, is the cardinality of a maximum $p$-independent set in $G$. In this paper, and motivated by the recent result that independence number is at most matching number for regular graphs~\cite{CaDaPe2020}, we investigate which values of the non-negative integers $p$, $q$, and $r$ have the property that $\alpha_p(G) \leq \alpha_q(L(G))$ for all r-regular graphs. Triples $(p, q, r)$ having this property are called \emph{valid $\alpha$-triples}.  Among the results we prove are:
\begin{itemize}
    \item $(p, q, r)$ is valid $\alpha$-triple for $p \geq 0$, $q \geq 3$ , and $r\geq 2$. 
    \item $(p, q, r)$ is valid $\alpha$-triple for $p \leq q < 3$ and $r\geq 2$. 
    \item $(p, q, r)$ is valid $\alpha$-triple for $p \geq 0$, $q = 2$, and $r$ even. 
    \item $(p, q, r)$ is valid $\alpha$-triple for $p \geq 0$, $q = 2$, and $r$ odd with $r = \max \Big \{ 3,  \frac{17(p+1)}{16}\Big \}$. 
\end{itemize}
We also show a close relation between undetermined possible valid $\alpha$-triples, the Linear Aboricity Conjecture, and the Path-Cover Conjecture. 
\end{abstract}


\section{Introduction}
Throughout this article, all graphs will be considered non-trivial, undirected, simple, and finite. In general, we will follow graph terminology and notation presented in~\cite{West}. Specifically, let $G$ be a graph with vertex set $V(G)$ and edge set $E(G)$. The order and size of $G$ will be denoted by $n(G) = |V(G)|$ and $m(G) = |E(G)|$, respectively. When the context is clear, we will use $n$ and $m$ instead of $n(G)$ and $m(G)$. We will denote the \emph{independence number} and the \emph{matching number} of $G$ by $\alpha(G)$ and $\mu(G)$, respectively. The independence and matching numbers are two of the oldest studied invariants in graph theory and are related by the following two theorems (among many others) given in~\cite{CaDaPe2020}.

\begin{theorem}[\cite{CaDaPe2020}]\label{thm:CaDaPe2020-general}
If $G$ is a graph with minimum degree $\delta$ and maximum degree $\Delta$, then
\[
\delta\alpha(G) \le \Delta\mu(G),
\]
and this bound is sharp. 
\end{theorem}

\begin{corollary}[\cite{CaDaPe2020}]\label{thm:CaDaPe2020-regular}
If $G$ is an $r$-regular graph with $r > 0$, then
\[
\alpha(G) \le \mu(G),
\]
and this bound is sharp. 
\end{corollary}
Notably, graphs attaining equality in Corollary~\ref{thm:CaDaPe2020-regular}, which was initially posed as a conjecture by the program \emph{TxGraffiti}~\cite{davila2024txgraffiti, TxGraffiti-website}, have been wholly characterized (see~\cite{cubic-characterization, regular-characterization}). The main aim of this paper is to expand and generalize Corollary~\ref{thm:CaDaPe2020-regular}.

\subsection{Generalized Independence and Matching}\label{subsection:generlized}
Let $p$ and $q$ denote non-negative integers, and let $r$ denote a positive integer. The \emph{line graph} of $G$, denoted $L(G)$, is the graph whose vertex set is the edge set of $G$, where two vertices in $L(G)$ are adjacent if and only if the edges they correspond to in $G$ are incident (share a vertex). A set of vertices $S \subseteq V(G)$ is called a \emph{$p$-independent set} of $G$ if each vertex in $S$ is adjacent to at most $p$ other vertices in $S$. The cardinality of a maximum $p$-independent set in $G$ is the \emph{$p$-independence number} of $G$, and is denoted by $\alpha_p(G)$. The $0$-independence number of $G$ is precisely the independence number; that is, $\alpha(G) = \alpha_0(G)$. A set of edges $M \subseteq E(G)$ is called a \emph{$q$-matching} of $G$ if each edge in $M$ is incident to at most $q$ other edges in $M$. The cardinality of a maximum $q$-matching of $G$ is called the \emph{$q$-matching number} of $G$, and is denoted by $\mu_q(G)$. The $0$-matching number of $G$ is precisely its matching number; that is, $\mu(G) = \mu_0(G)$.

It is clear from these definitions that the $q$-matching number of $G$ is identical to the $q$-independence number of its line graph $L(G)$. Namely, for every $q$, we have $\mu_q(G) = \alpha_q(L(G))$, and for this reason, we will now refer to the $q$-independence number of $L(G)$, without explicitly mentioning the $q$-matching number. Using this notation, we may rewrite the inequality stated in Corollary~\ref{thm:CaDaPe2020-regular} as $\alpha_0(G) \le \alpha_0(L(G))$, which implies $\alpha_p(G) \le \alpha_q(L(G))$ for any $r$-regular graph with $p = 0$, $q \geq 0$, and $r \geq 1$ -- a notion that corresponds to the validity of the triple $(p = 0, q \geq 0, r \geq 1)$ of integers which is generalized in the next section. This fact naturally leads to the following general question:
\begin{quest}\label{quest:main}
For which values of $p$, $q$, and $r$ is it true that all $r$-regular graphs satisfy the inequality
\[
\alpha_p(G) \le \alpha_q(L(G))?
\]
\end{quest}
When considering Question~\ref{quest:main}, it is only meaningful to consider $r \geq 1$, since the line graph $L(G)$ is not defined for any 0-regular graph $G$. 

\subsection{Contributions}\label{subsec:contributions}
Our first main contribution is to answer Question~\ref{quest:main} for many instances of $p$, $q$, and $r$. We call a triple $(p, q, r)$ that is valid for Question~\ref{quest:main} a \emph{valid $\alpha$-triple}; the only meaningful valid $\alpha$-triples being $(p, q, r \geq 1)$.  We prove a collection of valid $\alpha$-triples $(p, q, r \geq 1)$ and link well-known open problems in graph theory to the triples we could not prove.

As a consequence of considering Question~\ref{quest:main}, we give a natural and elegant generalization of Corollary~\ref{thm:CaDaPe2020-regular} as our second main contribution. We prove the following theorem, which, taking $p = 0$, reduces to the statement of Corollary~\ref{thm:CaDaPe2020-regular}.
\begin{theorem}\label{thm:generalized}
If $G$ is an $r$-regular graph with $r \geq 2$, then
\[
\alpha_p(G) \leq \alpha_p(L(G)).
\]
\end{theorem}

The remainder of this paper is organized as follows. In Section~\ref{sec:lemmas}, we provide known results, propositions, and lemmas needed to prove our main results. In Section~\ref{sec:main}, we answer Question~\ref{quest:main} for several varieties of triples. Section~\ref{sec:linear-arboricity} shows a relation between the unverified triples and the well-known Linear Arboricity and Path-Cover Conjectures. Finally, in Section~\ref{sec:conclusion}, we provide concluding remarks and pose several problems.

\section{Preliminaries}\label{sec:lemmas}
In this section, we give preliminary results needed for our study of Question~\ref{quest:main}.

\subsection{Graph Factor Lemmas}
Recall that a \emph{factor} of $G$ is a spanning subgraph $H$ of $G$. A \emph{$k$-factor} of $G$ is a spanning $k$-regular subgraph of $G$. The following two theorems are due to Petersen~\cite{Petersen} and Tutte~\cite{Tutte}, respectively.
\begin{thm}[\cite{Petersen}]\label{thm:petersen}
If $G$ is a $2r$-regular graph, then $G$ has a $2k$-factor for every integer $k$, where $0 < k < r$. 
\end{thm}

\begin{thm}[\cite{Tutte}]\label{thm:Tutte}
If $G$ is an $r$-regular graph and $k$ is an integer satisfying $0\leq k \leq  r$, then there exists a spanning subgraph $H$ of $G$ such that $k \leq d_H(v) \leq  k+1$ for each vertex $v \in V(G)$.
\end
{thm}

Next, recall that a \emph{$[k-1, k]$-factor} of $G$ is a spanning subgraph $H$ of $G$ such that for every vertex $v \in V(G)$, the degree of $v$ in $H$ satisfies $k-1 \le d_H(v) \le k$. A \emph{reduced $[k-1, k]$-factor}, denoted $R$-$[k-1, k]$-factor, is a spanning subgraph $H$ of $G$, such that for every vertex $v \in V(G)$, the degree of $v$ in $H$ satisfies $k-1 \le d_H(v) \le k$, and no two vertices of degree $k$ are adjacent in $H$. Note, every $[k-1, k]$-factor contains a $R$-$[k-1, k]$-factor, stated formally with the following result.

\begin{lem}\label{lem:factor}
Every $[k-1, k]$-factor contains a $R$-$[k-1, k]$-factor as a subgraph.
\end{lem}
\proof Let $G$ be a graph with a $[k-1, k]$-factor, say $H$. If no two vertices of degree $k$ in $H$ are adjacent in $H$, then we are done since this factor is already a reduced $[k-1, k]$-factor. If $H$ contains any adjacent vertices of degree $k$, choose two such vertices and delete the edge between them until no such pair remains. The subgraph obtained by removing these edges is spanning and an $R$-$[k-1,k]$ factor of $G$. \qed 

As a consequence of Theorem~\ref{thm:Tutte} and Lemma~\ref{lem:factor}, we have the following result. 
\begin{lem}\label{lem:kano}
If $G$ is an $r$-regular graph where $r \geq 3$ is an odd integer, then $G$ contains a $R$-$[2, 3]$-factor. 
\end{lem}

We will also need the following strengthening of Theorem~\ref{thm:Tutte} given by Kano~\cite{Kano}.
\begin{thm}[\cite{Kano}]\label{thm:Kano}
Let $j$ be a positive integer. If $k \leq \frac{2(2j+1)}{3}$, then every $(2j+1)$-regular graph has a $[k-1, k]$-factor each of whose components are regular. 
\end{thm}

The following result bounds the maximum degree of the line graph of a reduced $[k-1, k]$-factor and states a lower bound on the number of edges contained in such a factor. 
\begin{lem}\label{lem:factor-degree}
If $G$ is a graph of order $n$ and $H \subseteq G$ is a $R$-$[k-1, k]$-factor, then the following inequalities hold:
\begin{itemize}
    \item[(1)] $\Delta(L(H)) \leq 2k-3$
    \item[(2)] $|V(L(H))| = |E(H)| \geq \frac{n(k-1)}{2}$
\end{itemize} 
\end{lem}
\proof We first prove (1). Let $G$ be a graph and $H$ be a $R$-$[k-1, k]$-factor of $G$. Clearly, for every edge $e \in E(H)$ in the $R$-$[k-1, k]$-factor, the degree in the line graph is at most $k + (k-1) - 2 = 2k - 3$ since no two vertices of degree $k$ are adjacent in this factor. To see (2), note $\delta(H)  = k - 1$, which implies, 
\[
|E(H)| \geq \frac{n(k-1)}{2}.
\]
\qed

\subsection{Graph Matching Bounds}
In this section, we recall some well-known lower bounds relating to matching. First, recall the bound for 3-regular graphs given by Biedl et al. in~\cite{tight-matching-bound}.
\begin{thm}[\cite{tight-matching-bound}]\label{thm:tight-matching-bound}
If $G$ is a 3-regular graph of order $n$, then 
\[
\alpha_0(L(G)) = \mu(G) \geq \frac{4n - 1}{9},
\]
and this bound is sharp. 
\end{thm}

Theorem~\ref{thm:tight-matching-bound} is a special instance of the more general bound given by Henning and Yeo in~\cite{matching-lower-bounds}.
\begin{thm}[\cite{matching-lower-bounds}]\label{thm:henning-yeo}
For odd $r$, if $G$ is a connected $r$-regular graph of order $n$, then
\[
\alpha_0(L(G)) = \mu(G) \geq \frac{(r^3 - r^2 - 2)n - 2r + 2}{2(r^3 - 3r)},
\]
and this bound is sharp. 
\end{thm}

\subsection{Graph Independence Bounds}
We will also need to use lower and upper bounds on the generalized independence number we present in this section. To begin, we first prove a lemma that can also be derived from the \emph{degree sequence index strategy} (\emph{DSI-strategy}) introduced in~\cite{DSI}; see also~\cite{hansberg-pepper} for another alternative proof and an extensive study of generalized independence in graphs. However, we provide direct proof from first principles to stay as self-contained as possible.
\begin{lem}\label{lem:general-upper}
If $G$ is an $r$-regular graph of order $n$ and $0 \le p \le r$, then
\begin{equation}\label{eq:1}
  \alpha_p(G) \le \frac{nr}{2r-p}, 
\end{equation}
and this bound is sharp. If in addition $p < r$, then 
\begin{equation}\label{eq:2}
  \alpha_p(G) \le \frac{nr}{r+1}.  
\end{equation}
\end{lem}

\proof If $r = p$, then (\ref{eq:1}) is trivial, as $\alpha_p(G) \leq n$. Thus, we may assume $r > p$. Let $G$ be an $r$-regular graph of order $n$ with $r = p + t$, where $t \ge 1$. Next, let $A \subseteq V(G)$ be a maximum $p$-independent set of $G$ and let $B = V(G) \setminus A$. Thus, $|A| = \alpha_p(G)$ and $|B| = n - \alpha_p(G)$. Let $e(A, B)$ denote the set of edges with one endpoint in $A$ and the other in $B$. Since $A$ is a $p$-independent set, if $v \in A$, then $v$ can be incident with at most $p$ edges whose other endpoint is also in $A$. Thus, each vertex in $A$ contributes at least $r-p = t$ edges to the count $e(A, B)$. Hence, $e(A, B) \ge t \cdot |A|$. On the other hand, since $G$ is $r$-regular, if $v \in B$, then $v$ can be incident with at most $r$ edges whose other endpoint is in $A$. Thus, each vertex in $B$ contributes at most $r$ edges to the count $e(A, B)$, and so $e(A, B) \le r \cdot |B|$. Hence,
\[
t \cdot |A| = t \cdot (n - |B|) \le e(A, B) \le r \cdot |B|.
\]
Rearranging this expression, we obtain
\[
\alpha_p(G) = |A| \le \frac{nr}{r+t} = \frac{nr}{2r-p},
\]
establishing inequality (\ref{eq:1}). Using $p < r$ establishes inequality (\ref{eq:2}).

We next show that inequality (\ref{eq:1}) is sharp. Consider the bipartite graph $G = [A, B]$ with $|A| = 2m$, where each vertex in $A$ has degree $r-p$, and each vertex in $B$ has degree $r$. Next, pack a $p$-factor into $A$ and denote the resulting graph by $G'$. Note that for $|A| = 2m$, the complete graph $K_{2m}$ has a 1-factorization, and so, you take $p$ 1-factors to pack into $A$. From the cardinalities of $A$ and $B$, it follows that $(r-p) \cdot |A| = r \cdot |B|$, while $|A| + |B| = n$. Thus,
\[
(r-p) \cdot |A| = r \cdot (n - |A|).
\]
Rearranging this equality, we obtain
\[
|A| = \frac{nr}{2r-p}.
\]
Since $A$ is a $p$-independent set of $G'$, inequality (\ref{eq:1}) is sharp. \qed

\begin{remark}\label{remark}
A bipartite graph $G$, as claimed in the proof above, is constructed, for example, by taking $2t$ copies of $K_{r,r - p}$ and can be made connected, if one wishes, by switching operations. We mention here for later use that if $2r-p$ is not divisible by three, then for $2t$, where $t$ is not divisible by three, the number of vertices of the graph $G$ constructed above is not divisible by three.
\end{remark}

We next focus on bounding the generalized independence number of the line graph $L(G)$ in terms of the order of the graph $G$. To start, the following lemma bounds the 1-independence number of the line graph. 
\begin{lem}\label{lem:alpha_1-trivial}
If $G$ is a graph of order $n$, then 
\[
\alpha_1(L(G)) \leq \frac{2n}{3},
\]
and this bound is sharp.
\end{lem}
\proof Observe that a maximum 1-independent set in $L(G)$ corresponds to a subgraph of the form $x P_3 \cup yP_2$ in $G$ which gives a 1-independent set in $L(G)$ of cardinality $2x+y$. Clearly, $3x +2y \leq n$. Hence, 
\[
\alpha_1(L(G)) = 2x + y \leq 2x + \frac{n-3x}{2} = \frac{n + x}{2} \leq \frac{2n}{3},
\]
as $x \leq \frac{n}{3}$. To see that this bound is sharp, take $G$ to be any graph with a $P_3$-packing. \qed

Next we focus on bounding the the 2-independence number of the line graph $L(G)$ of $G$ from below in terms of the order $n$ of $G$.
\begin{lem}\label{lem:3-parts}
Let $G$ be an $r$-regular graph with $r \geq 2$ and order $n$. 
\begin{enumerate}
    \item[(1)] If $r$ is even, then $\alpha_2(L(G)) = n$. 
    \item[(2)] If $r \geq 3$, then $\alpha_2(L(G)) \geq \frac{(17n - 2)}{18}$.
\end{enumerate}
\end{lem}
\proof 
Let $G$ be an $r$-regular graph with $r \geq 2$ and order $n$. We first prove (1). Let $r$ be even. Then, by Theorem~\ref{thm:petersen}, $G$ has a 2-factor, say $F$, with exactly $n$ edges in $G$. The edges of $F$ in $G$ form a 2-independent set in $L(G)$, which implies $\alpha_2(L(G)) \geq n$.  Suppose $\alpha_2(L(G)) \geq n + 1$ for some $r$-regular graph $G$. In $G$, there exists a set $F$ of at least $n + 1$ edges, such that the corresponding vertices in $L(G)$, denoted by $L(F)$, form a set of $n + 1$ vertices in $L(G)$ that induces a subgraph with maximum degree at most two. This is impossible as in the subgraph induced by $F$ in $G$ there must be a cycle with an attached edge forming degree three in $L(F)$. Thus, $\alpha_2(L(G)) = n$ completing the proof of (1). 


We next prove (2). Let $r = 3$ and observe that in this case the graph $G$ has $m(G) = \frac{3}{2}n$ edges. We next remove edges from $G$ so that the resulting graph contains no vertex of degree three. By Theorem~\ref{thm:tight-matching-bound}, there exists a matching, say $M \subseteq E(G)$, such that $|M| \geq \frac{4n-1}{9}$. Note that the number of $M$-saturated vertices in $G$ is $2|M|$. Therefore, there is at most $n - 2|M|$ vertices which are not the endpoint of any edge in $M$. For each vertex in $G$ which is not the endpoint of an edge in $M$, we delete one edge. Hence, we delete at most
\[
|M| + (n - 2|M|) = n - |M| \leq n - \frac{4n-1}{9} = \frac{5n + 1}{9},
\]
edges before we arrive at a graph with no degree three vertex. That is, the resulting subgraph of $G$, say $H$, satisfies $\Delta(H) = 2$, and also, 
\[
m(H) \geq \frac{3}{2}n - \frac{5n + 1}{9} = \frac{27n-10n-2}{18} = \frac{17n-2}{18}.
\]
Since $\Delta(H) = 2$, the edges of $H$ form a 2-independent set in $L(G)$. Thus, 
\[
\alpha_2(L(G)) \geq m(H) \geq \frac{17n-2}{18}.
\]

We now consider $r \geq 5$. By Theorem~\ref{thm:Kano} we may assume $G$ contains a $[2,3]$-factor that splits $V(G) = V(G_2) \cup V(G_3)$ where $G_2$ is 2-regular graph on $n_2$ vertices and $G_3$ is 3-regular graph on $n_3$ vertices and $n_2 +  n_3 = n$. Thus,
\begin{align*}
\alpha_2(L(G)) &\geq \alpha_2(L(G_2)) + \alpha_2(L(G_3)) \\
               &\geq n_2  +  \frac{17n_3 - 2}{18}  \\
               &= \frac{ 18(n - n_3) + (17n_3  - 2)}{18} \\
               &=  \frac{18n - n_3 - 2}{18} \\
               &\geq \frac{17n - 2}{18},
\end{align*}
which completes the proof of statement (2). \qed 

\section{Main Results}\label{sec:main}
This section presents our main results addressing Question~\ref{quest:main}. To establish a foundation for the general case, we begin by examining the almost trivial scenario of 1-regular graphs.
\begin{prop}\label{prop:r=1}
If $G$ is a $1$-regular graph, then
\[
\alpha_p(G) \leq \alpha_q(L(G)),
\]
with equality if and only if $p = 0$.
\end{prop}
\proof For $r = 1$, the 1-regular graph $G$ of order $n$, consists of $n/2$ disjoint copies of $K_2$, and so, the line graph $L(G)$ consists of $n/2$ disjoint copies of $K_1$. Therefore, for $r = 1$, $\alpha_0(G) = \alpha_q(L(G))$ for $q \geq 0$. However, for $p \geq 1$ we observe $\alpha_p(G) = n$. Thus, $\alpha_p(G) \leq \alpha_q(L(G))$ holds for 1-regular graphs if and only if $p = 0$. \qed

\begin{cor}\label{cor:r=1}
$(p, q \geq 0, r = 1)$ is a valid $\alpha$-triple if and only if $p = 0$.
\end{cor}

The intriguing cases of Question~\ref{quest:main} arise for $r$-regular graphs with $r \geq 2$, and the remainder of this section is dedicated to addressing these cases. Each subsection concludes with a corollary or observation on valid or invalid $\alpha$-triples, and we highlight open cases where Question~\ref{quest:main} remains unresolved.


\subsection{When $r$ is Even}
When $r \geq 2$ is even, we may apply Petersen's Theorem and obtain our first nontrivial valid $\alpha$-triple. 
\begin{thm}\label{thm:even}
If $G$ is an $r$-regular graph with $r$ even, then,
\[
\alpha_p(G) \leq \alpha_q(L(G)),
\]
whenever $p \geq 0$ and $q \geq 2$. 
\end{thm}
\proof Let $G$ be an $r$-regular graph with $r$ chosen to be even. In this case, Theorem~\ref{thm:petersen} implies that $G$ has a 2-factor, say $H \subseteq G$, and moreover, $|E(H)| = n$. Since $H$ is a 2-factor in $G$, $L(H)$ is a 2-regular graph in $L(G)$. Thus, $V(L(H))$ is a 2-independent set, implying $\alpha_p(G) \leq n = \alpha_2(L(G))$ for all $p \geq 0$, and we are done. \qed 

\begin{cor}\label{cor:r-is-even}
$(p \geq 0 , q \geq 2, r = 2k)$ is a valid $\alpha$-triple for all integers $k \geq 1$.
\end{cor}

By Corollary~\ref{cor:r-is-even}, it follows that Question~\ref{quest:main} is resolved whenever $r \geq 2$ is even and $q\geq 2$. The remaining cases when $r \geq 2$ is even are $q = 0$ and when $q = 1$, of which we address in general in the next two subsections.

\subsection{When $q = 0$}
For $q= 0$ it is clear that $\alpha_0(L(G)) = \mu(G) \leq n/2$. Thus, for fixed $r$ and $p \geq r$, there are no valid $\alpha$-triples since $\alpha_p(G) = \alpha_r(G) = n > \alpha_0(L(G))$ in this case.  Next recall, by Lemma~\ref{lem:general-upper}, that for $0 \leq p \leq  r$, 
\[
\alpha_p(G) \leq  \frac{nr}{2r - p},
\]
and that this bound is sharp -- constructions achieving sharpness are given in the proof of Lemma~\ref{lem:general-upper}. Take $G$ to be one of the graphs achieving equality in Lemma~\ref{lem:general-upper}, and observe that whenever $p > 0$,
\[
\alpha_p(G) = \frac{nr}{2r - p } > \frac{n}{2} \geq \alpha_0(L(G)) = \mu(G).
\]
Hence, $(p \geq 1, q = 0, r \geq 1)$ is never a valid $\alpha$-triple. Combined with Corollary~\ref{thm:CaDaPe2020-regular}, this implies that the only valid $\alpha$-triples with $q = 0$ are $(p = 0, q = 0, r \geq 1)$.
\begin{cor}\label{cor:q=0}
$(p = 0 , q = 0, r \geq 1)$ is a valid $\alpha$-triple.
\end{cor}


\subsection{When $q \geq 3$}
With the following theorem we show all triples of the form $(p, q \geq 3, r \geq 2)$ are valid $\alpha$-triples. 
\begin{thm}\label{lem:q>=3}
If $G$ is an $r$-regular graph with $r \geq 2$ and $q \geq 3$, then
\[
\alpha_p(G) \leq \alpha_q(L(G)),
\]
for $p \geq 0$.
\end{thm}
\proof Let $G$ be an $r$-regular graph of order $n$ with $r \geq 2$. 
If $r$ is even, then we are done by Theorem~\ref{thm:even}. Thus, we may assume that $r$ is odd. Since $\alpha_p(G) \leq n$ for all $p$, and $\alpha_3(L(G)) \leq \alpha_q(L(G))$ for all $q \geq 3$ by monotonicity, it suffices to show that $\alpha_3(L(G)) \geq n$. 

If $r = 3$, $G$ is a $[2, 3]$-factor (with no degree 2 vertices). Then, by Lemma~\ref{lem:factor}, $G$ has a $R$-$[2, 3]$-factor. If $r \geq 5$, then Lemma~\ref{lem:kano} implies that $G$ has a $R$-$[2, 3]$-factor. In both cases, we are assured that $G$ contains a $R$-$[2, 3]$-factor, say $H \subseteq G$. Then, by Lemma~\ref{lem:factor-degree}, $|E(H)|\geq n$ and $\Delta(L(H)) \leq 3$. Thus, $V(L(H))$ forms a 3-independent set in $L(G)$ with cardinality at least $n$. That is,
\[
n \leq |V(L(H))| = |E(H)| \leq \alpha_3(L(G)).
\]
and the proof is complete. \qed

\begin{cor}\label{cor:q>=3}
$(p \geq 0 , q \geq 3, r \geq 2)$ is a valid $\alpha$-triple.
\end{cor}

\subsection{When $p \leq q < 3$}
We now consider Question~\ref{quest:main} for the difficult case when $p \leq q < 3$. More specifically, the following theorem shows all triples of the form $(p \leq q, q < 3, r \geq 2)$ are valid $\alpha$-triples.  
\begin{thm}\label{thm:main}
If $G$ is an $r$-regular graph with $r \geq 2$, then
\[
\alpha_p(G) \le \alpha_q(L(G)),
\]
whenever $p \leq q < 3$. 
\end{thm}

\proof Since the case of $p=q=0$ was already established, we only address the cases with $q \geq 1$ in this proof. Thus we aim to establish that the inequality is true when $(p,q) \in \{ (0,1),(0,2), (1,1), (1,2), (2,2) \}$. Furthermore, we will frequently appeal to the monotonicity property of the $k$-independence number -- namely, that $\alpha_k(G) \leq \alpha_{k+1}(G)$ for every non-negative integer $k$. To this end, assume $G$ is an $r$-regular graph with $r \ge 2$. We consider two separate cases: $r = 2$ and $r \geq 3$.

\medskip
\noindent \textbf{Case 1:} $r = 2$.

In this case, $G$ is a $2$-regular graph, and therefore a disjoint union of cycles. It follows immediately that $L(G) = G$, and thus we have
\[
\alpha_2(G) = |V(G)| = \alpha_2(L(G)).
\]
Consequently, by the monotonicity of the $k$-independence number,
\[
\alpha_0(G) \leq \alpha_1(G) \leq \alpha_2(G) = \alpha_2(L(G)),
\]
implying that $(p,q) \in \{(0,2),(1,2),(2,2)\}$ are valid pairs. Additionally, we note that $\alpha_1(G) = \alpha_1(L(G))$, so 
\[
\alpha_0(G) \leq \alpha_1(G) = \alpha_1(L(G)),
\]
validating the pairs $(p,q) \in \{(0,1),(1,1)\}$, again by appeal to monotonicity of $\alpha_k$.

\medskip
\noindent \textbf{Case 2:} $r \ge 3$.

Here, we only need to explicitly address the cases $(p,q) = (1,1)$ and $(p,q) = (2,2)$; as the remaining pairs follow by monotonicity of $k$-independence.

\smallskip
\noindent \textbf{Subcase 2.1:} $p = q = 1$.

Let $A$ be a largest $1$-independent set in $G$, chosen to maximize the number of edges in the induced subgraph $G[A]$. Thus, no additional edge can be added to $A$ without violating $1$-independence, and clearly, $|A| = \alpha_1(G)$. Let $B = V(G) \setminus A$. By definition of $1$-independence, each component of $A$ is either an isolated edge ($K_2$) or an isolated vertex ($K_1$).


Since $r\ge 3$, every vertex in $A$ has at least one neighbor in $B$, ensuring $B$ is nonempty. The choice of $A$ yields two useful properties:

\begin{itemize}
    \item \textbf{Property I:} Each $K_2$ component contributes exactly $2(r-1)$ edges to $B$, while each $K_1$ component contributes exactly $r$ edges to $B$.
    \item \textbf{Property II:} Each vertex in $B$ has at most $r$ neighbors in $A$.
\end{itemize}

Next, let $S$ be the set of components of $A$ and construct a bipartite multigraph $H$ with vertex set $S \cup B$, where each vertex in $S$ corresponds uniquely to one component of $A$. A vertex $u \in S$ is connected to $v \in B$ by edges whose multiplicity equals the number of vertices in the corresponding component of $A$ adjacent to $v$.


Let $D \subseteq S$ and consider the edges in the induced bipartite subgraph $H[D, N(D)]$. Property~{II} implies there are at most $r|N(D)|$ edges going from $N(D)$ to $D$. However, Property~{I} guarantees at least $r|D|$ edges going from $D$ to $N(D)$, giving the inequality
\[
r|N(D)| \ge r|D|.
\]
This implies $|D| \leq |N(D)|$. Thus, Hall’s theorem ensures the existence of a matching from $S$ into $B$, assigning to each component of $A$ a distinct neighbor in $B$. Using this matching, each component of type $K_2$ extends to $P_3$, and each component of type $K_1$ extends to $P_2$, forming a $1$-independent set of the same cardinality as $A$ in $L(G)$. Thus, the pair $(1,1)$ is valid and by monotonicity, $(0,1)$ also follows.

\smallskip
\noindent \textbf{Subcase 2.2:} $p = q = 2$.

Let $A$ be a largest $2$-independent set in $G$, again chosen to maximize the number of edges in $G[A]$. Thus, no further edge can be included without breaking $2$-independence, and $|A|=\alpha_2(G)$. Set $B=V(G)\setminus A$. By maximality, the induced subgraph on $A$ consists of disjoint cycles, paths, and isolated vertices. Because $r \ge 3$, each vertex in $A$ must have at least one neighbor in $B$, ensuring $B$ is nonempty. Since cycle components remain unchanged under the line-graph transformation (i.e., $L(C_k)=C_k$), we may restrict attention solely to non-cycle components. This choice of $A$ yields two useful properties:

\begin{itemize}
    \item \textbf{Property III:} The two end-points of each path component collectively contribute exactly $2(r-1)$ edges to $B$, and each isolated vertex ($K_1$) contributes exactly $r$ edges to $B$.
    \item \textbf{Property IV:} Each vertex in $B$ has at most $r$ neighbors in $A$.
\end{itemize}

Let $S$ be the set of non-cycle components of $A$ and construct a bipartite multigraph $H$ with vertex set $S \cup B$, where vertices in $S$ correspond uniquely to the non-cycle components of $A$. Connect $u \in S$ to $v \in B$ with multiplicity equal to the number of endpoints (or single vertex, in the $K_1$ case) adjacent to $v$ in $G$.


Let $D\subseteq S$ and consider the edges in the induced subgraph $H[D,N(D)]$. Property~{IV} implies there are at most $r|N(D)|$ edges going from $N(D)$ to $D$. However, Property~{III} ensures at least $r|D|$ edges going from $D$ to $N(D)$. This gives the inequality
\[
r|N(D)| \ge r|D|,
\]
which implies $|D| \leq |N(D)|$. Thus, Hall’s theorem ensures the existence of a matching from $S$ into $B$, assigning to each non-cycle component of $A$ a distinct neighbor in $B$. This matching enables each path component to extend by one additional vertex and each isolated vertex to form an edge, which when taken together with the cycle components of $A$, produces a $2$-independent set with the same cardinality as $A$ in $L(G)$. Hence, the pair $(2,2)$ is valid, and monotonicity immediately validates the pairs $(1,2)$ and $(0,2)$.

\medskip
All cases and subcases are now complete, establishing the theorem.
\qed

\begin{cor}\label{cor:p<=q<3}
$(p \leq q, q < 3, r \geq 2)$ is a valid $\alpha$-triple. 
\end{cor}

By Corollary~\ref{cor:q>=3} and Corollary~\ref{cor:p<=q<3}, we deduce the appealing corollary below, which confirms Corollary~\ref{thm:generalized}. Namely, for all $r$-regular graphs with $r \geq 2$, the $p$-independence number is at most the $p$-matching number for all non-negative integers $p$.
\begin{cor}[Corollary~\ref{thm:generalized}]
If $G$ is an $r$-regular graph, with $r \geq 2$, and $p$ is a non-egative integer, then 
\[
\alpha_p(G) \leq \mu_p(G) = \alpha_p(L(G)).
\]
\end{cor}

\subsection{When $q = 2$}
This section addresses some of the remaining cases not covered by Theorem~\ref{thm:main} when $q=2$, namely when $p \geq 3$.

\begin{thm}\label{17p}
If $G$ is an $r$-regular graph with odd $r$ such that $r \geq \max \Big \{ 3, \frac{17(p+1)}{16} \Big \}$, then
\[
\alpha_p(G) \le \alpha_q(L(G)),
\]
whenever $p \geq 0$ and $q \geq 2$. 
\end{thm}
\proof Let $G$ be an $r$-regular graph of order $n$ with $ r \geq 3$, where $r$ is odd. If $q \geq 3$, then we are done by Theorem~\ref{lem:q>=3}. If $q = 2$ and $p \in \{0, 1, 2 \}$, then we are also done by Theorem~\ref{thm:main}. Thus,  we may assume $q = 2$ and $p \geq 3$. Hence, $r$ is odd with $r \geq  17(p+1)/16 > 3$. Next observe, 

\[
\begin{array}{rcl}
\frac{rn}{2r-p} & \leq & \frac{17n-2}{18} \\
& \Updownarrow & \\
18rn & \leq & (2r-p)(17n-2) \\
& \Updownarrow & \\
18rn & \leq & 34rn - 4r - 17pn + 2p \\
& \Updownarrow & \\
16rn & \geq & 17pn - 2p + 4r.
\end{array}
\]

Thus, since $r \geq  17(p+1)/16$, 
\[
16rn \geq 17n(p+1) \geq 17pn - 2p + 4r \iff  17np + 17n  \geq 17pn - 2p + 4r.
\]
Hence, $17n \geq 4r – 2p$, which is a true inequality since $n \geq r +1$ (recalling $G$ is $r$-regular). By Lemma~\ref{lem:general-upper} and Lemma~\ref{lem:3-parts}, we finally observe, 
\[
\alpha_p(G) \leq \frac{rn}{2r-p} \leq \frac{17n-2}{18} \leq \alpha_2(L(G)), 
\]
completing the proof of the theorem. \qed 

\begin{cor}\label{cor:p>=3}
If $r$ is odd, then $(p \geq 0 , q = 2, r \geq \max \Big \{ 3, \frac{17(p+1)}{16} \Big \})$ is a valid $\alpha$-triple.
\end{cor}

\subsection{When $q = 1$}
For $q = 1$, many non-valid $\alpha$-triples exist. For example, consider any $r$-regular graph, say $G$, which attains the equality
\[
\alpha_p(G) = \frac{rn}{2r - p},
\]
where we refer the reader to the proof of Lemma~\ref{lem:general-upper} for one such construction of the graph $G$. Next recall Lemma~\ref{lem:alpha_1-trivial}, which states that 
\[
\alpha_1(L(G)) \leq \frac{2n}{3},
\]
for any graph $G$ of order $n$. We observe,
\[
\frac{rn}{2r - p} > \frac{2n}{3} \iff 3rn > 4rn -2pn \iff 3r > 4r -2p \iff 2p > r,
\]
Thus, choosing $2p > r$ yields
\[
\alpha_1(L(G)) < \frac{rn}{2r - p} = \alpha_p(G).
\]
Therefore, no triple of the form $(p \geq 0 , q = 1, 1 \leq r < 2p)$ is a valid $\alpha$-triple. Furthermore, by Remark~\ref{remark}, we also observe that for $2r-p$ not divisible by 3, the same argument above gives that $( p \geq 0 , q = 1, 1 \leq r \leq 2p)$ is also not a valid $\alpha$-triple.

\section{Relation to Linear Arboricity and Path-Cover Conjectures}\label{sec:linear-arboricity}
The only open cases for $q = 2$ are when $p \geq 3$ and odd $r$ with $3 \leq  r <  17(p+1)/6$. As we shall see, the missing cases are closely related to two famous conjectures: The Linear Arboricity Conjecture and the Path-cover conjecture.  

To illustrate this, recall that in 1970, Harary introduced the concept of \emph{linear arboricity} as a covering invariant on graphs~\cite{Harary1}. A \emph{linear forest} is a graph in which each component is a path, and the linear arboricity of a graph $G$, denoted $la(G)$, is the minimum number of linear forests whose union is $G$. The following conjecture~\cite{Harary2}, open since 1980, is of fundamental importance in research on linear arboricity.

\begin{conj}[The Linear Arboricity Conjecture~\cite{Harary2}]\label{Linear-Arboricity-Conjecture} 
The linear arboricity of an $r$-regular graph is $\lceil \frac{r+1}{2} \rceil$.
\end{conj}

Assuming the Linear Arboricity Conjecture is true, if $G$ is an $r$-regular graph with odd $r$ and $r \ge 3$, then $G$ can be covered by $\frac{r+1}{2}$ linear forests. This covering contains, by averaging, a linear forest $F$ with at least,
\[
\frac{|E(G)|}{la(G)} = \frac{(rn/2)}{((r+1)/2)} = \frac{rn}{r+1}
\]
edges, which form a 2-independent set in $L(G)$. Therefore, if the Linear Arboricity Conjecture is true, then $\alpha_q(L(G)) \ge \alpha_2(L(G)) \ge \frac{rn}{r+1}$ for all $q \ge 2$. Moreover, if $r > p$, then Lemma~\ref{lem:general-upper} (2) implies $\alpha_p(G) \le \frac{rn}{r+1}$. These two facts lead to the following observation.

\begin{cor}
 If $G$ is an $r$-regular graph with $r \geq  p+1 \geq 3$ and we assume the Linear Arboricity Conjecture true, then 
 \[
 \alpha_p(G) \leq \alpha_q(L(G)),
 \]
 for all $q \geq 2$.
\end{cor}

The Linear Arboricity Conjecture has been confirmed for $r \in \{1, 2, 3, 4, 5, 6, 8, 10\}$ (see~\cite{Guldan}). Assuming the linear arboricity, the triple $( p \geq 0 ,q = 2 , r \geq p+1)$ is a valid $\alpha$-triple, while $(p \geq 0, q = 2, 3 \leq r \leq p)$ is not a valid $\alpha$-triple because for $p \geq r \geq 1$, $\alpha_p(G) = n$ while $\alpha_2(L(G)) = n$ if and only if $G$ has 2-factor. As is well known (see~\cite{Boll1985}), there exist $r$-regular graphs with odd $r \geq 3$ that have no $k$-factor for $1 \leq k \leq r - 1$, and in particular, no 2-factor. For such graphs, $\alpha_2(L(G)) < n$. For example, let $k$ be odd and construct a graph $G$ of order $k+2$ whose complement has degree sequence $2, 1, ...., 1$. Next take $k$ copies of $G$ together with an isolated vertex, say $u$, where we attach $u$ to vertices of $G$ with degree $k-1$ in $G$; a graph that has no nontrivial factor. Thus, If the Linear Arboricity Conjecture is true for all odd $r$, then the determination of valid $\alpha$-triples in case $q = 2$ is complete.




The \emph{path cover number} of $G$, denoted by $\rho(G)$, is the minimum number of vertex-disjoint paths required to cover the vertices of $G$. In~\cite{Magnant}, Magnant and Martin conjectured the following, a conjecture which remains open in general but has been proven for $r$-regular graphs with $0 \le r \le 6$ (see~\cite{r=6, Magnant, approximatepathdecompositionsregular}).

\begin{conj}[The Path-Cover Conjecture~\cite{Magnant}]
If $G$ is an $r$-regular graph of order $n$, then $\rho(G) \le \frac{n}{r+1}$.
\end{conj}

Suppose $G$ is a graph of order $n$ that is coverable by $k$ paths, say $P_{1}, \dots, P_{k}$, each with order $n_i$, for $i = 1, \dots, k$, respectively. Since every vertex in $G$ is in exactly one path in this covering, we have $n_1 + \dots + n_k = n$. Furthermore, since the path $P_{j}$ has exactly $n_j-1$ edges, the total number of edges in this path covering is $(n_1-1) + \dots + (n_k -1) = n - k$. Hence, $\alpha_2(L(G)) \ge n - \rho(G)$. Thus,
if we assume the Path-Cover Conjecture is true and $G$ is an $r$-regular graph with $r > p$, then 
\[
\alpha_2(L(G)) \ge n - \frac{n}{r+1} = \frac{rn}{r+1} \ge \alpha_2(G),
\]
where the right-hand side of this inequality follows by Lemma~\ref{lem:general-upper} (2).

Assuming the Path-Cover Conjecture true, $(p \geq 0 ,  q = 2 ,  r \geq  p+1)$ is a valid $\alpha$-triple, while $( p \geq 0, q = 2  , 3 \leq  r \leq p)$ is not a valid $\alpha$-triple as for $p \geq  r \geq \alpha_p(G) = n$ while $\alpha_2(L(G)) = n$ if and only if $G$ has 2-factor. But again, and as before, it is well known that there are $r$-regular graphs with odd $r \geq 3$ that have no $k$-factor for $1 \leq k \leq r - 1$ and for these graphs $\alpha_2(L(G)) < n$. Thus, if the Path-Cover Conjecture is true for all odd $r$, then the determination of valid $\alpha$-triples in case $q = 2$ is complete. 


\section{Concluding Remarks}\label{sec:conclusion}
In this section we summarize our main results concerning Question 1 and offer several problems; Table~\ref{tab:example} below gives all answers presented in this paper. 

\begin{table}[h!]
\centering
\begin{tabularx}{\textwidth}{|>{\centering\arraybackslash}m{2cm}|>{\centering\arraybackslash}m{2cm}|>{\centering\arraybackslash}m{2cm}|>{\centering\arraybackslash}m{2cm}|>{\centering\arraybackslash}X|}
\hline
\textbf{$q$} & \textbf{$r$} & \textbf{$p$} & \textbf{Valid $\alpha$-triple?} & \textbf{Comments} \\
\hline
$q \geq 3$ & $r \geq 2$ & $p \geq 0$ & Yes & Theorem~\ref{lem:q>=3} \\
\hline
$q = 2$ & $r$ even & $p \geq 0$ & Yes &  Theorem~\ref{thm:even}\\
\hline
$q \leq 2$ & $r \geq 2$ & $p \leq q$ & Yes &  Theorem~\ref{thm:main} \\
\hline
$q = 2$ & $r \geq 3$ & $p \geq 3$ & Yes* &  *for $r$ odd and $r \geq \frac{17(p+1)}{16}$, Theorem~\ref{17p}\\
\hline
$q = 2$ & $r \geq 3$ & $p \geq 3$ & No &  $r$ odd and $r \leq p$\\
\hline
$q = 2$ & $r \geq 3$ & $p \geq 3$ & Unknown & $r$ odd and $p + 1 \leq r < \frac{17(p+1)}{16}$. Yes, if Linear Arboricity or Path-Cover Conjecture \\
\hline 
$q \geq 1$ & $r = 1$ & $p = 0$ & Yes & $\alpha_p(G)=\alpha_q(L(G))=\frac{n}{2}$ \\
\hline 
$q \geq 1$ & $r = 1$ & $p \geq 1$ & No & $\alpha_p(G)=n$ while $\alpha_q(L(G))<n$ \\
\hline
$q = 1$ & $r = 2$ & $p \geq 2$ & No &  $\alpha_p(G)=n$ while $\alpha_q(L(G))<n$\\
\hline
$q = 1$ & $r \geq 3$ & $2p > r$ & No &  In the case that $2r-p$ is not divisible by 3, $2p\geq r$\\
\hline
$q = 0$ & $r \geq 1$ & $p \geq 1$ & No & Section 3.2 \\
\hline
$q = 0$ & $r \geq 1$ & $p = 0$ & Yes &  Corollary~\ref{thm:CaDaPe2020-general} (see also~\cite{CaDaPe2020})\\
\hline

\end{tabularx}
\caption{Summary of triples considered.}
\label{tab:example}
\end{table}

The results summarized in Table~\ref{tab:example} raise the following two problems.
\begin{problem}\label{prob:1}
Is is true that $(p \geq 3, q = 2, p+1 \leq r < \frac{17(p+1)}{16})$ are all valid $\alpha$-triples?
\end{problem}

\begin{problem}\label{prop:2}
For $q = 1$ and $p \geq 2$ and $r \geq 3$ determine which $(p \geq 2, q = 1, r \geq 3)$ are valid $\alpha$-triples. 
\end{problem}

In another direction we propose the following: In~\cite{choi}, it is proved that every 3-regular graph contains a 2-regular subgraph of cardinality at least $\min\{n, \frac{5(n+2)}{6}\}$. For our purpose to get lower-bound on $\alpha_2(L(G))$, where $G$ is an $r$-regular graph with $r \geq 3$, we considered the case when $G$ had no 2-factor, and used a lower bound on the number of edges in a subgraph $H$ with $\Delta(H) \leq 2$, and proved the lower bound $\frac{17n-2}{18}$. 

Recall for a given graph $G$ and $k \geq 0$, $\mu_k(G)$ is the size of a largest subset $H$ of edges with the property that each edge in $H$ is incident to at most $k$ other edges in $H$, so that clearly, $\alpha_k(L(G))=\mu_k(G)$. If we define $\mu_{r,k}(n)$ as the minimum value of $\mu_k(G)$ over all $r$-regular graphs $G$ of order $n$, then we can state another problem for future research.

\begin{problem}
For a given $r \geq 3$, determine, or give upper or lower approximations, for $\mu_{r,k}(n)$. 
\end{problem}
 
 
\section*{Acknowledgments}
The authors would like to thank the anonymous referees for their careful reading of the
manuscript and for their insightful comments and suggestions.


\medskip
\begin{thebibliography}{99}

\bibitem{Harary2} J. Akiyama, G. Exoo, and F. Harary, Covering and packing in graphs III. Cyclic and acyclic invariants. \emph{Math. Slovaca}, \textbf{30} (1980), 405--417.

\bibitem{tight-matching-bound} T. Biedl, E. D. Demaine, C. A. Duncan, R. Fleischer, and S. G. Kobourov, Tight bounds on maximal and maximum matchings, \emph{Discrete Math.}, \textbf{285(1-3)} (2004), 7--15.

\bibitem{Boll1985} B. Bollobas, Akira Saito, and N.C. Wormald, Regular factors of regular graphs, \emph{J. of Graph Theory}, \textbf{9(1)} (1985), 97--103.

\bibitem{Caro-k-domination} Y. Caro, On the $k$-domination and $k$-transversal numbers of graphs and hypergraphs, \emph{Ars Combin.}, \textbf{29 C} (1990), 49--55.

\bibitem{CaDaPe2020} Y. Caro, R. Davila, and R. Pepper, New results relating independence and matchings, \textit{Discuss. Math. Graph Theory}, \textbf{42(3)} (2022), 921--935.

\bibitem{DSI} Y. Caro and R. Pepper, Degree sequence index strategy, \emph{Australas. J. Combin.}, \textbf{59(1)} (2014), 1--23.

\bibitem{choi} I. Choi, R. Kim, A. V. Kostochka, B. Park, and D. B. West, Largest 2-Regular Subgraphs in 3-Regular Graphs, \emph{Graphs Combin.}, \textbf{35} (2019), 805--813. 


\bibitem{davila2024txgraffiti} R. Davila, Automated conjecturing in mathematics with \emph{TxGraffiti}, arXiv preprint arXiv:2409.19379, 2024. 


\bibitem{TxGraffiti-website} R. Davila, \emph{TxGraffiti} Interactive Website. Available at: \url{https://txgraffiti.streamlit.app}, 2024. Accessed: 2024-07-23.

\bibitem{Guldan} F. Guldan, The linear arboricity of 10-regular graphs, \emph{Math. Slovaca}, \textbf{36(3)} (1986), 225--228.

\bibitem{prob-k-dom} D. Rautenbach and L. Volkmann, New bounds on the $k$-domination number and the $k$-tuple domination number, \emph{Appl. Math. Lett.}, \textbf{20(1)} (2007), 98--102.

\bibitem{r=6} U. Feige and E. Fuchs, On the path partition number of 6-regular graphs, \emph{J. Graph Theory}, \textbf{101(3)} (2022), 345--378. 

\bibitem{hansberg-pepper} A. Hansberg and R. Pepper, On $k$-domination and $j$-independence in graphs, \emph{Discrete Appl. Math.}, \textbf{161} (2013), 1472--1480.

\bibitem{Harary1} F. Harary, Covering and packing in graphs {I}, \emph{Ann. N. Y. Acad. Sci.}, \textbf{175} (1970), 198--205. 

\bibitem{matching-lower-bounds} M. A. Henning and A. Yeo, Tight lower bounds on the size of a maximum matching in a regular graph, \emph{Graphs Combin.}, \textbf{23} (2007), 647--657. 


\bibitem{Magnant} C. Magnant and D.M. Martin, A note on the path cover number of regular graphs, \emph{Australas. J. Combin.}, \textbf{43} (2009), 211--217.

\bibitem{cubic-characterization} E. Mohr and D. Rautenbach, Cubic graphs with equal independence number and matching number, \emph{Discrete Math.}, \textbf{344(1)} (2021), 112178.

\bibitem{approximatepathdecompositionsregular} R. Montgomery, A. Müyesser, A. Pokrovskiy, and B. Sudakov, Approximate path decompositions of regular graphs, arXiv preprint arXiv:2406.02514, 2024. 

\bibitem{Kano} M. Kano, Factors of regular graphs, \emph{J. Comb. Theory Series B}, \textbf{41(1)} (1986), 27--36.

\bibitem{Petersen} J. Petersen, Die Theorie der regularen Graph, \emph{Acta Math.}, \textbf{15} (1891), 193--220.

\bibitem{Plummer2} M. Plummer, \emph{Factors and Factorization}. 403-430. Handbook of Graph Theory ed. J. L. Gross and J. Yellen. CRC Press, 2003, ISBN: 1-58488-092-2.

\bibitem{Tutte} W.T. Tutte, The subgraph problem, 
\emph{Annals of Discrete Math.}, \textbf{3} (1978), 289--295.

\bibitem{West} D. B. West, Introduction to Graph Theory, Second Edition, Prentice-Hall (20010. ISBN: 0-13-014400-2 (print)

\bibitem{regular-characterization} Z. Yang and H. Lu, Regular graphs with equal matching number and independence number, \emph{Discrete Appl. Math.}, \textbf{310} (2022), 86--90.

\end{thebibliography}
\end{document}